\documentclass[12pt,a4paper,reqno]{amsart}
\usepackage{amsmath}
\usepackage{amsfonts}
\usepackage{amssymb}
\usepackage{bm}

\newcommand\R{{\mathbf{R}}}

\newcommand\Z{{\mathbf{Z}}}

\newcommand\bigO{{\mathcal{O}}}

\theoremstyle{plain}
  \newtheorem{theorem}[subsection]{Theorem}

\theoremstyle{remark}
  \newtheorem{remark}[subsection]{Remark}

\theoremstyle{definition}

\include{psfig}

\begin{document}

\title[Global regularity for logarithmically supercritical Navier-Stokes]{Global regularity for a logarithmically supercritical hyperdissipative Navier-Stokes equation}
\author{Terence Tao}
\address{Department of Mathematics, UCLA, Los Angeles CA 90095-1555}
\email{tao@math.ucla.edu}

\vspace{-0.3in}
\begin{abstract}
Let $d \geq 3$.  We consider the global Cauchy problem for the generalised Navier-Stokes system
\begin{align*}
\partial_t u + (u \cdot \nabla) u &= - D^2 u - \nabla p \\
\nabla \cdot u &= 0 \\
u(0,x) &= u_0(x)
\end{align*}
for $u: \R^+ \times \R^d \to \R^d$ and $p: \R^+ \times \R^d \to \R$, where $u_0: \R^d \to \R^d$ is smooth and divergence free, and $D$ is a Fourier multiplier whose symbol $m: \R^d \to \R^+$ is non-negative; the case $m(\xi) = |\xi|$ is essentially Navier-Stokes.  It is folklore (see e.g. \cite{kp}) that one has global regularity in the critical and subcritical hyperdissipation regimes $m(\xi) = |\xi|^\alpha$ for $\alpha \geq \frac{d+2}{4}$.  We improve this slightly by establishing global regularity under the slightly weaker condition that $m(\xi) \geq |\xi|^{(d+2)/4}/g(|\xi|)$ for all sufficiently large $\xi$ and some non-decreasing function $g: \R^+ \to \R^+$ such that $\int_1^\infty \frac{ds}{sg(s)^4} = +\infty$. 
In particular, the results apply for the logarithmically supercritical dissipation $m(\xi) := |\xi|^{\frac{d+2}{4}} / \log(2 + |\xi|^2)^{1/4}$.
\end{abstract}

\maketitle

\section{Introduction}

Let $d \geq 3$.  This note is concerned with solutions to the generalised Navier-Stokes system
\begin{equation}\label{ns}
\begin{split}
\partial_t u + (u \cdot \nabla) u &= -D^2 u - \nabla p \\
\nabla \cdot u &= 0 \\
u(0,x) &= u_0(x)
\end{split}
\end{equation}
where $u: \R^+ \times \R^d \to \R^d$, $p: \R^+ \times \R^d \to \R$ are smooth, and $u_0: \R^d \to \R^d$ is smooth, compactly supported, and divergence-free, and $D$ is a Fourier multiplier\footnote{The exact definition of the Fourier transform is not relevant for this discussion, but for sake of concreteness one can take $\hat f(\xi) := \int_{\R^d} f(x) e^{-i x \cdot \xi}\ dx$.} whose symbol $m: \R^d \to \R^+$ is non-negative; the case $m(\xi) = |\xi|$ is essentially the Navier-Stokes system, while the case $m=0$ is the Euler system.

For $d \geq 3$, the global regularity of the Navier-Stokes system is of course a notoriously difficult unsolved problem, due in large part to the supercritical nature of the equation with respect to the energy $E(u(t)) := \int_{\R^d} |u(t,x)|^2\ dx$.  This supercriticality can be avoided by strengthening the dissipative symbol $m(\xi)$, for instance setting $m(\xi) := |\xi|^\alpha$ for some $\alpha > 1$.  This \emph{hyper-dissipative} variant of the Navier-Stokes equation becomes subcritical for $\alpha > \frac{d+2}{4}$ (and critical for $\alpha = \frac{d+2}{4}$), and it is known that global regularity can be recovered in these cases; see \cite{kp} for further discussion.  For $1 \leq \alpha < \frac{d+2}{4}$, only partial regularity results are known; see \cite{ckn} for the $\alpha=1$ case and \cite{kp} for the $\alpha > 1$ case.

The purpose of this note is to extend the global regularity result very slightly into the supercritical regime:

\begin{theorem}\label{main}  Suppose that $m$ obeys the lower bound
\begin{equation}\label{mangle}
 m(\xi) \geq |\xi|^{(d+2)/4} / g(|\xi|)
\end{equation}
for all sufficiently large $|\xi|$, where $g: \R^+ \to \R^+$ is a non-decreasing function such that
\begin{equation}\label{gss}
 \int_1^\infty \frac{ds}{s g(s)^4} = \infty.
\end{equation}
Then for any smooth, compactly supported initial data $u_0$, one has a global smooth solution to \eqref{ns}.
\end{theorem}

Note that the hypotheses are for instance satisfied when 
\begin{equation}\label{mxi}
m(\xi) := |\xi|^{(d+2)/4} / \log^{1/4}(2+|\xi|^2),
\end{equation}
thus
$$ |D|^2 = \frac{(-\Delta)^{(d+2)/4}}{\log^{1/2}(2 - \Delta)}.$$
Analogous ``barely supercritical'' global regularity results were established for the nonlinear wave equation recently in \cite{tao-supercrit}, \cite{roy1}, \cite{roy2}.

The argument is in fact quite simple, being based on the classical energy method and Sobolev embedding.  The basic point is that whereas in the critical and subcritical cases one can get an energy inequality of the form
$$ \partial_t \| u(t) \|_{H^k(\R^d)}^2 \leq C a(t) \| u(t) \|_{H^k(\R^d)}^2$$
for some locally integrable function $a(t)$ of time, a constant $C$, and some large $k$, which by Gronwall's inequality is sufficient to establish a suitable \emph{a priori} bound, in the logarithmically supercritical case \eqref{mxi} one instead obtains the slightly weaker inequality
$$ \partial_t \| u(t) \|_{H^k(\R^d)}^2 \leq C a(t) \| u(t) \|_{H^k(\R^d)}^2 \log(2 + \| u(t) \|_{H^k(\R^d)})$$
(thanks to an endpoint version of Sobolev embedding, closely related to an inequality of Br\'ezis and Wainger \cite{bw}), which is still sufficient to obtain an \emph{a priori} bound (though one which is now double-exponential rather than single-exponential, cf. \cite{bkm}).

\begin{remark} It may well be that the condition \eqref{gss} can be relaxed further by a more sophisticated argument.  Indeed, the following heuristic suggests that one should be able to weaken \eqref{gss} to $\int_1^\infty \frac{ds}{s g(s)^2} = \infty$, thus allowing one to increase the $1/4$ exponent in \eqref{mxi} to $1/2$.  Consider a blowup scenario in which the solution blows up at some finite time $T_*$, and is concentrated on a ball of radius $1/N(t)$ for times $0 < t < T_*$, where $N(t) \to \infty$ as $t \to T_*$.  As the energy of the fluid must stay bounded, we obtain the heuristic bound $u(t) = O(N(t)^{d/2})$ for times $0 < t < T_*$.  In particular, we expect the fluid to propagate at speeds $O( N(t)^{d/2} )$, leading to the heuristic ODE $\frac{d}{dt} \frac{1}{N(t)} = O( N(t)^{d/2} )$ for the radius $1/N(t)$ of the fluid.  Solving this ODE, we are led to a heuristic upper bound $N(t) = O( (T_*-t)^{\frac{2}{d+2}})$ on the blowup rate.  On the other hand, from the energy inequality
$$ 2 \int_0^{T_*} \int_{\R^d} |Du(t,x)|^2\ dx dt \leq \int_{\R^d} |u_0(x)|^2\ dx$$
one is led to the heuristic bound
$$ \int_0^{T_*} \frac{1}{N(t)^{(d+2)/2} g(N(t))^2}\ dt < \infty.$$
This is incompatible with the upper bound $N(t) = O( (T_*-t)^{\frac{2}{d+2}})$ if $\int_1^{T_*} \frac{ds}{s g(s)^2} = \infty$.  Unfortunately the author was not able to make this argument precise, as there appear to be multiple and inequivalent ways to rigorously define an analogue of the ``frequency scale'' $N(t)$, and all attempts of the author to equate different versions of these analogues lost one or more powers of $g(s)$.

To go beyond the barrier $\int_1^\infty \frac{ds}{s g(s)^2} = \infty$ (with the aim of getting closer to the Navier-Stokes regime, in which $g(s) = s^{1/4}$ in three dimensions), the above heuristic analysis suggests that one would need to force the energy to not concentrate into small balls, but instead to exhibit turbulent behaviour.
\end{remark}

\subsection{Acknowledgements}

The author is supported by NSF Research Award DMS-0649473, the NSF Waterman award and a grant from the MacArthur Foundation.  

\section{Proof of theorem}

We now prove the theorem.  Let $k$ be a large integer (e.g. $k := 100d$ will suffice).  

Standard energy method arguments (see e.g. \cite{kato}) show that if the initial data is smooth and compactly supported, then either a smooth $H^\infty$ solution exists for all time, or there exists a smooth solution up to some blowup time $0 < T_* < \infty$, and $\|u(t)\|_{H^k(\R^d)} \to \infty$ as $t \to T_*$.  Thus, to establish global regularity, it suffices to prove an \emph{a priori bound} of the form
$$ \| u(t) \|_{H^k(\R^d)} \leq C( k, d, \| u_0 \|_{H^k(\R^d)}, T, g )$$
for all $0 \leq t \leq T < \infty$ and all smooth $H^\infty$ solutions $u: [0,T] \times \R^d \to \R^d$ to \eqref{ns}, where $C( k, d, \| u_0 \|_{H^k(\R^d)}, T, g )$ is a quantity depending only on $k$, $d$, $\| u_0 \|_{H^k(\R^d)}$, $T$, and $g$.  

We now fix $u_0, u, T$, and let $C$ denote any constant depending on $k$, $d$, $\| u_0 \|_{H^k(\R^d)}$, $T$, and $g$ (whose value can vary from line to line).  Multiplying the Navier-Stokes equation by $u$ and integrating by parts, we obtain the well-known energy identity
$$ \partial_t \int_{\R^d} |u(t,x)|^2\ dx = - 2 a(t)$$
where 
\begin{equation}\label{adu}
a(t) := \| Du \|_{L^2(\R^d)}^2
\end{equation}
(note that the pressure term $\nabla p$ disappears thanks to the divergence free condition); integrating this in time, we obtain the energy dissipation bound
\begin{equation}\label{energy-dissipate}
\int_0^T a(t)\ dt \leq C.
\end{equation}

Now, we consider the higher energy 
\begin{equation}\label{ek-def}
E_k(t) := \sum_{j=0}^k \int_{\R^d} |\nabla^j u(t,x)|^2\ dx.
\end{equation}
Differentiating \eqref{ek-def} in time and integrating by parts, we obtain
$$
\partial_t E_k(t) = - 2 \sum_{j=0}^k \| \nabla^j D u(t) \|_{L^2(\R^d)}^2 - 2 \sum_{j=0}^k \int_{\R^d} \nabla^j u(t,x) \cdot \nabla^j ( (u \cdot\nabla) u)(t,x)\ dx;$$
again, the pressure term disappears thanks to the divergence-free condition.  For brevity we shall now drop explicit mention of the $t$ and $x$ variables.

We apply the Leibniz rule to $\nabla^j( (u \cdot \nabla) u)$.  There is one term involving $j+1^{th}$ derivatives of $u$, but the contribution of that term vanishes by integration by parts and the divergence free property.  The remaining terms give contributions of the form
$$ \sum_{j=0}^k \sum_{1 \leq j_1,j_2 \leq j: j_1+j_2=j+1} \int_{\R^d} \bigO( \nabla^j u \nabla^{j_1} u \nabla^{j_2} u )\ dx$$
where $\bigO( \nabla^j u \nabla^{j_1} u \nabla^{j_2} u )$ denotes some constant-coefficient trilinear combination of the components of $\nabla^j u$, $\nabla^{j_1} u$, and $\nabla^{j_2} u$ whose explicit form is easily computed, but is not of importance to our argument.  We can integrate by parts using $D$ and $D^{-1}$ and then use Cauchy-Schwarz to bound
$$
\int_{\R^d} \bigO( \nabla^j u \nabla^{j_1} u \nabla^{j_2} u )\ dx \leq \| (1+D) \nabla^j u \|_{L^2(\R^d)} \| (1+D)^{-1}( \bigO( \nabla^{j_1} u \nabla^{j_2} u ) ) \|_{L^2(\R^d)}.$$
By the arithmetic mean-geometric mean inequality we then have
$$
\int_{\R^d} \bigO( \nabla^j u \nabla^{j_1} u \nabla^{j_2} u )\ dx \leq c \| (1+D) \nabla^j u \|_{L^2(\R^d)}^2 + \frac{1}{c} \| (1+D)^{-1}( \bigO( \nabla^{j_1} u \nabla^{j_2} u ) ) \|_{L^2(\R^d)}^2$$
for any $c>0$.  Finally, from the triangle inequality, \eqref{ek-def}, and the fact that $D$ commutes with $\nabla^j$, we have
$$ \| (1+D) \nabla^j u \|_{L^2(\R^d)}^2 \leq C ( \| \nabla^j D u \|_{L^2(\R^d)}^2 + E_k ).$$
Putting this all together and choosing $c$ small enough, we conclude that
\begin{equation}\label{chock}
 \partial_t E_k \leq C E_k + C \sum_{1 \leq j_1 \leq j_2 \leq k: j_1+j_2\leq k+1} 
\| (1+D)^{-1} ( \bigO( \nabla^{j_1} u \nabla^{j_2} u ) ) \|_{L^2(\R^d)}^2.
\end{equation}
To estimate this expression, we introduce a parameter $N > 1$ (depending on $t$) to be optimised later, and divide $(1+D)^{-1} = (1+D)^{-1} P_{\leq N} + (1+D)^{-1} P_{> N}$, where $P_{\leq N}$ and $P_{>N}$ are the Fourier projections to the regions $\{ \xi: |\xi| \leq N \}$ and $\{ \xi: |\xi| > N \}$.  

We first deal with the low-frequency contribution to \eqref{chock}.  From Plancherel's theorem and \eqref{mangle} we see that
$$ \| (1+D)^{-1} P_{\leq N} ( \bigO( \nabla^{j_1} u \nabla^{j_2} u ) ) \|_{L^2(\R^d)} \leq C g(N) \| \langle \nabla \rangle^{-(d+2)/4} \bigO( \nabla^{j_1} u \nabla^{j_2} u ) \|_{L^2(\R^d)} $$
where $\langle \nabla \rangle^{-(d+2)/4}$ is the Fourier multiplier with symbol $\langle \xi \rangle^{-(d+2)/4}$, where $\langle \xi\rangle := (1+|\xi|^2)^{1/2}$.  Applying Sobolev embedding, we can bound the right-hand side by
$$ \leq C g(N) \| |\nabla^{j_1} u| |\nabla^{j_2} u| \|_{L^{4d/(3d+2)}(\R^d)}.$$
By H\"older's inequality and the Gagliardo-Nirenberg inequality, we can bound this by
$$ \leq C g(N) \| \nabla u \|_{L^{4d/(d+2)}(\R^d)} \| \nabla^{j_1+j_2-1} u \|_{L^2(\R^d)} $$
which by \eqref{ek-def} is bounded by
$$ \leq C g(N) \| \nabla u \|_{L^{4d/(d+2)}(\R^d)} E_k^{1/2}.$$
Next, we partition
$$ \| \nabla u \|_{L^{4d/(d+2)}(\R^d)} \leq \| \nabla P_{\leq N} u \|_{L^{4d/(d+2)}(\R^d)} + \| \nabla P_{> N} u \|_{L^{4d/(d+2)}(\R^d)}.$$
From Sobolev embedding and Plancherel, \eqref{mangle} and \eqref{adu}, we have
\begin{align*}
\| \nabla P_{\leq N} u \|_{L^{4d/(d+2)}(\R^d)} &\leq C \| \langle \nabla \rangle^{(d+2)/4} P_{\leq N} u \|_{L^2(\R^d)} \\
&\leq C g(N) (1+a(t))^{1/2}.
\end{align*}
Meanwhile, from Sobolev embedding we have
$$ \| \nabla P_{> N} u \|_{L^{4d/(d+2)}(\R^d)} \leq \frac{1}{N} E_k^{1/2}$$
(say) if $k$ is large enough.  Putting this all together, we see that the low-frequency contribution to \eqref{chock} is
$$ \leq C g(N)^2 E_k [ g(N)^2 (1+a(t)) + \frac{1}{N^2} E_k ].$$

Next, we turn to the high-frequency contribution to \eqref{chock}.  From Plancherel, H\"older's inequality, and \eqref{ek-def} we have
\begin{align*}
\| (1+D)^{-1} P_{\geq N} ( \bigO( \nabla^{j_1} u \nabla^{j_2} u ) ) \|_{L^2(\R^d)}
&\leq C g(N) N^{-(d+2)/4} \| |\nabla^{j_1} u| |\nabla^{j_2} u| \|_{L^2(\R^d)} \\
&\leq C g(N) N^{-(d+2)/4} \| \nabla^{j_1} u\|_{L^\infty(\R^d)} E_k^{1/2},
\end{align*}
while from Sobolev embedding and \eqref{ek-def} we see (for $k$ large enough) that
$$ \| \nabla^{j_1} u\|_{L^\infty(\R^d)} \leq C E_k^{1/2}.$$
Thus the high-frequency contribution to \eqref{chock} is
$$ \leq C g(N)^2 N^{-(d+2)/2} E_k^2.$$
Putting this all together, we conclude that
$$ \partial_t E_k \leq C g(N)^2 E_k [ g(N)^2 (1+a(t)) + \frac{1}{N} E_k ].$$
We now optimise in $N$, setting $N := 1 + E_k$, to obtain
$$ \partial_t E_k \leq C g(1+E_k)^4 E_k (1+a(t)).$$
From \eqref{energy-dissipate}, \eqref{gss} and separation of variables we see that the ODE
$$ \partial_t E = C g(1+E)^4 E (1+a(t))$$
with initial data $E(0) \geq 0$ does not blow up in time.  Also, from \eqref{ek-def} we have $E_k(0) \leq C$.  A standard ODE comparison (or continuity) argument then shows that $E_k(t) \leq C(T)$ for all $0 \leq t \leq T$, and the claim follows.

\begin{remark} It should be clear to the experts that the domain $\R^d$ here could be replaced by any other sufficiently smooth domain, e.g. the torus $\R^d/\Z^d$, using standard substitutes for the Littlewood-Paley type operators $P_{\leq N}$, $P_{>N}$ (e.g. one could use spectral projections of the Laplacian).  We omit the details.
\end{remark}

\end{document}